\documentclass[12pt]{article}
\usepackage{amsmath}
\usepackage{CJK}
\usepackage{amsfonts}
\usepackage{mathrsfs}
\usepackage{amssymb}
\usepackage[usenames]{color}
\def\RED{}\def\REMARK#1{}
\def\g{\mathfrak{g}}

\def\l{\lambda}
\def\m{\mu}

\def\p{\partial}

\def\Vir{{\rm Vir}}
\def\sc{\scriptstyle}
\def\ssc{\scriptscriptstyle}

\def\cl{\centerline}

\def\vs{\vspace*}

\def\ni{\noindent}
\def\ptl{\partial}

\def\Z{\mathbb{Z}{\ssc\,}}

\def\C{\mathbb{C}{\ssc\,}}

\numberwithin{equation}{section}
\newtheorem{theo}{Theorem}[section]
\newtheorem{defi}[theo]{Definition}
\newtheorem{rema}[theo]{Remark}
\newtheorem{lemm}[theo]{Lemma}

\def\fg{{\mathfrak{g}}}
\begin{document}
\begin{CJK*}{GBK}{song}
\setcounter{page}{0}
%
%

\cl{{\large \bf Filtered Lie conformal algebras whose associated
graded algebras }}\cl{{\large\bf  are isomorphic to that of general
conformal algebra $gc_1$} \footnote {Supported by NSF grant
no.~10825101 and 11001200 of China and the Program for
Young Excellent Talents in Tongji University.\\\indent
Correspending author: X. Yue (xiaoqingyue@tongji.edu.cn)}}
\vs{6pt}

\cl{Yucai Su$^{\,\dag,\,\ddag}$,  \ \ Xiaoqing Yue$^{\,\dag}$} \cl{\footnotesize $^{\dag\,}$Department of Mathematics, Tongji University, Shanghai
200092, China}
\cl{\footnotesize $^{\ddag\,}$Wu Wen-Tsun Key Laboratory of \vs{-2pt}Mathematics}
\cl{\footnotesize   University of Science and
Technology of China, Hefei 230026, China}
\cl{\footnotesize E-mail: ycsu@tongji.edu.cn, xiaoqingyue@tongji.edu.cn}
\vs{6pt} \vs{10pt}\par

{\small \parskip .005 truein \baselineskip 3pt \lineskip 3pt

\noindent{{\bf Abstract.} Let $G$ be a filtered Lie conformal algebra whose associated
graded conformal algebra is isomorphic to that of general conformal
algebra $gc_1$. In this paper, we prove that $G\cong gc_1$ or ${\rm gr\,}gc_1$ (the associated
graded conformal algebra of $gc_1$), by making use of some results on the second cohomology groups
of the conformal algebra $\fg$ with coefficients in its  module $M_{b,0}$ of rank $1$, where $\fg=\Vir{\sc\!}\ltimes{\sc\!} M_{a,0}$ is the semi-direct sum of the Virasoro conformal algebra $\Vir$ with its module $M_{a,0}$. Furthermore, we prove that
${\rm gr\,}gc_1$ does not have a nontrivial representation on a finite $\C[\partial]$-module, this provides an
example of a finitely freely generated simple Lie conformal algebra of linear
growth that cannot be embedded into the general conformal algebra $gc_N$ for any $N$.
 \vs{5pt}

\noindent{\bf Key words:} Lie conformal algebras, filtered algebras, general conformal algebras, second cohomology groups}

\noindent{\it Mathematics Subject Classification (2000):} 17B10,
17B65,  17B68.}
\parskip .001 truein\baselineskip 8pt \lineskip 8pt

\vs{6pt}
\par

\cl{\bf\S1. \
Introduction}\setcounter{section}{1}\setcounter{equation}{0}\setcounter{theo}{0}

The notion of  conformal algebras, introduced  in [\ref{K1}],
encodes an axiomatic description of the operator product expansion
(or rather its Fourier transform) of chiral fields in conformal
field theory. Conformal algebras play important roles in quantum
field theory and vertex operator algebras, which is also an adequate
tool for the study of infinite-dimensional Lie algebras and
associative algebras (and their representations), satisfying the
locality property in [\ref{K2}]. There has been a great deal of work
towards understanding of the algebraic structure underlying the
notion of the operator product expansion of chiral fields of a
conformal field theory (e.g., [\ref{BPZ}, \ref{K3}]). The singular
part of the operator product expansion encodes the commutation
relations of fields, which leads to the notion of  Lie conformal
algebras.

The structure theory, representation theory and cohomology theory of
finite Lie conformal algebras has been developed in the past few
years (e.g., [\ref{BKLR}, \ref{BKV}, \ref{CK}--\ref{FKR}, 
\ref{X1}--\ref{Z2}
]). Simple finite Lie conformal algebras were classified in
[\ref{DK}], which shows that a simple finite conformal algebra is
isomorphic either to the Virasoro conformal algebra or to the
current Lie conformal algebra Cur$\,\g$ associated to a simple
finite-dimensional Lie algebra $\g$. The theory of conformal modules
has been developed in [\ref{CK}] and of their extensions in
[\ref{CKW}]. In particular, all finite simple irreducible
representations of simple finite Lie conformal algebras were
constructed in [\ref{CK}]. The  cohomology theory of conformal
algebras with coefficients in an arbitrary module has been developed
in [\ref{BKV}]. However, the structure theory, representation theory
and cohomology theory of simple infinite Lie conformal algebras is
 far from being well developed.

In order to better understand the 
theory of simple infinite Lie conformal
algebras, it is very natural
 to first study some important examples. As is pointed out in [\ref{K1}],
 the general Lie conformal algebras and the associative conformal algebras
 are the most important examples of simple infinite conformal algebras. It
 is well-known that the general Lie conformal algebra $gc_N$ plays the same
  important role in the theory of Lie conformal algebras as the general Lie
   algebra $gl_N$ does in the theory of Lie algebras: any module $M=\C[\p]^{N}$
    over a Lie conformal algebra $R$ is obtained via a homomorphism $R\rightarrow gc_N$. Thus
the study of Lie conformal algebras $gc_N$ has drawn some authors'
attentions (e.g., [\ref{BKL1}, \ref{BKL2}, \ref{S1}]).

In order to study simple infinite Lie conformal algebras, it is also
very important to find the filtered infinite Lie conformal algebras
whose associated graded Lie conformal algebras are some known
infinite Lie conformal algebras. The determination of filtered
algebras with associated graded algebras isomorphic to some known
algebras is in general a highly nontrivial problem, as can be seen from
examples in [\ref{Kac-f1}, \ref{Kac-f2}].
 Thus one of our motivations in the present paper
is to investigate some simple infinite Lie conformal algebras by
determining filtered Lie conformal algebras whose associated graded
conformal algebras are isomorphic to that of general conformal
algebra $gc_N$. Due to the reason stated in Remark \ref{rema-gc-N}(2), we have to treat
$gc_1$ separately. Thus in the present paper, we shall only consider the case $gc_1$.
The consideration of  the general case $gc_N$ for $N>1$ will be our next goal.

The main result in the present paper is the following theorem\vs{-5pt}.
\begin{theo}{\bf(Main Theorem\vspace*{-3pt})}\label{theo}\begin{enumerate}\parskip-3pt\item[\rm(1)]
Let $G$ be a filtered Lie conformal algebra whose associated graded
conformal algebra is isomorphic to that of general conformal algebra
$gc_1$. Then   $G\cong gc_1$ or ${\rm gr\,}gc_1$ $($the associated graded
conformal algebra of $\,gc_1)$.
\item[\rm(2)] The graded conformal algebra ${\rm gr\,}gc_1$ does not have a nontrivial representation on any finite $\C[\partial]$-module.
In particular, ${\rm gr\,}gc_1\not\cong gc_1$, and ${\rm gr\,}gc_1$ is a finitely freely generated simple Lie conformal algebra of linear growth that is not embedded into $gc_N$ for any $N$.
\end{enumerate}
\end{theo}

Theorem \ref{theo}(1) will follow from some computations (see Section 4) and some results (Theorem \ref{Prop1}) on the second cohomology groups of the conformal algebra $\fg$ with coefficients in its  modules $M_{b,0}$ of rank $1$, where $\fg=\Vir\ltimes M_{a,0}$ is the semi-direct sum of the Virasoro conformal algebra $\Vir$ with its module $M_{a,0}$. Here in general, $\Vir$ and $M_{\Delta,\alpha}$ (which is a simple $\Vir$-module if and only if $\Delta\ne0$) are defined \vs{-5pt}by
\begin{eqnarray}
\label{Vir}&&
\Vir=\C[\ptl]L:\ \ \ \ L_\l L=(\ptl+2\l)L.\\
\label{M-a-b}&&
M_{\Delta,\alpha}=\C[\ptl]v:\ \ L_\l v=(\alpha+\ptl+\Delta\l)v.
\end{eqnarray}

We will prove Theorem \ref{theo}(2) in Section 5.
\vs{5pt}

\cl{\bf\S2. \ The general Lie conformal algebra \vspace*{-7pt}$gc_1$
}\setcounter{section}{2}\setcounter{equation}{0} \setcounter{theo}{0}

\begin{defi}\rm A
{\it Lie conformal algebra } is a $\C[\p]$-module $A$ with a
$\l$-bracket $[x_{\l}y]$ which defines a linear map $A\times A
\rightarrow A[\l]$, where $A[\l]=\C[\l]\otimes A$ is the space of
polynomials of $\l$ with coefficients in $A$, such that for \vs{-5pt}$x,y,z\in A$,\RED{
\begin{eqnarray}&\!\!\!\!\!\!\!\!\!\!\!\!\!\!\!\!\!\!\!\!\!\!\!\!&
[\p x_{\l}y]=-\l [x_\l y],\ \ \ \  [x_\l \p y]=(\p+\l)[x_\l y] \ \ \   \mbox{(conformal  sesquilinearity)},\\
&\!\!\!\!\!\!\!\!\!\!\!\!\!\!\!\!\!\!\!\!\!\!\!\!&\label{J_a_b_}
[x_{\l}y]=-[y_{-\l-\p} x] \ \ \   \mbox{(skew-symmetry)},\\
&\!\!\!\!\!\!\!\!\!\!\!\!\!\!\!\!\!\!\!\!\!\!\!\!&\label{J_a_b_c}
[x_{\l}[y_\m z]]=[[x_\l y]_{\l+\m} z]+[y_\m[x_\l z]] \ \ \
\mbox{(Jacobi identity)}.
\end{eqnarray}
}
\end{defi}

The {\it general Lie conformal algebra} $gc_N$ can be defined as the
infinite rank
$\C[\p]$-module  $\C[\ptl,x]\otimes gl_N$,
 with the $\l$-bracket
\begin{equation}\label{gc-N}
[{f(\ptl,x)A}_{\,\l\,} g(\ptl,x)B]{\sc\!}={\sc\!}
f(-\l,x{\sc\!}+{\sc\!}\ptl{\sc\!}+{\sc\!}\l)
g(\ptl{\sc\!}+{\sc\!}\l,x)AB{\sc\!}-{\sc\!}f(-\l,x)g(\ptl{\sc\!}
+{\sc\!}\l,x{\sc\!}-{\sc\!}\l)BA,
\end{equation}
for $f(\ptl,x),g(\ptl,x)\in\C[\ptl,x],$ $A,B\in gl_N$, where $gl_N$ is the space of $N\times N$ matrices, and we have identified $f(\ptl,x)\otimes A$ with $f(\ptl,x)A$.
If we set $J_A^n=x^nA$, \vs{-10pt}then
\begin{equation*}
[{J_{A}^{m}}_{\l} J_B^{n}]=\mbox{$\sum\limits_{s=0}^{m}$}\binom{m}{s}
(\l+\p)^{s}J_{AB}^{m+n-s}-\mbox{$\sum\limits_{s=0}^{n}$} \binom{n}{s}
(-\l)^{s}J_{BA}^{m+n-s}\vs{-5pt},
\end{equation*}
for $m,n\in \Z_+,$ $A,B\in gl_N,$ where
$\big(^m_{\,s}\big)=m(m-1)\cdot\cdot\cdot(m-s+1)/{s}!$ if $s\geq 0$
and $\big(^m_{\,s}\big)=0$ otherwise, is the binomial coefficient.
The formal
distribution Lie algebra corresponding to $gc_N$ is the well-known
Lie algebra $\mathcal{D}^N$ of $N\times N$-matrix differential
operators on the circle.

In particular,  the {\it general Lie conformal algebra} $gc_1$ is the
 infinite rank free  $\C[\p]$-module $\C[\ptl,x]$
with a generating set
$\{J_{n}=x^{n+1}\mid -1\le n\in\Z\}$, such \vs{-5pt}that
\begin{equation}\label{gc1-bracket-2}
[J_{{m}_{\l}} J_{n}]=\mbox{$\sum\limits_{s=0}^{m}$}\binom{m+1}{s+1}
(\l+\p)^{s+1}J_{m+n-s}-\mbox{$\sum\limits_{s=0}^{n}$}\binom{n+1}{s+1}
(-\l)^{s+1}J_{m+n-s}\vs{-5pt},
\end{equation}for $-1\le m,n\in \Z$.
Naturally, $gc_1$ is a filtered algebra with filtration
\begin{eqnarray}\label{gc-1-fil}
\{0\}{\sc\!}={\sc\!}gc_1^{(-2)}{\sc\!}\subset
 {\sc\!}gc_1^{(-1)}{\sc\!}\subset{\sc\!}\cdots{\sc\!}\subset{\sc\!} gc_1\mbox{ with }
gc_1^{(n)}{\sc\!}={\sc\!}\mbox{span}\{J_i\,|\,-1{\sc\!}\leq{\sc\!} i{\sc\!}\leq{\sc\!} n\}
\mbox{ for } n{\sc\!}\geq{\sc\!}-1,
\end{eqnarray}
such that
\begin{equation}\label{gc1-bracket-3}
[J_{i_{\l}} J_j]\equiv \Big((i+j+2)\l+(i+1)\p\Big) J_{i+j}
\mbox{~~(~mod~}gc_1^{(i+j-1)})\mbox{ \ for \ }i,j\ge-1.
\end{equation}

\begin{defi}\rm\label{defi-gr}
Let $G=\cup_{i\in \Z} G^{(i)}$ be a Lie conformal algebra with  a filtration $$
\cdots\subseteq G^{(-1)}\subseteq G^{(0)}\subseteq G^{(1)}\subseteq \cdots,$$ such that
$[G^{(i)}{}_{_{\l}}G^{(j)}]{\sc\!}\subseteq{\sc\!} G^{(i+j)}$ for $i,j{\sc\!}\in{\sc\!} \Z$. Denote ${\rm gr\,}G{\sc\!}={\sc\!}\oplus _{i\in\Z}\overline{G^{(i)}}$, where $\overline{G^{(i)}}=G^{(i)}/G^{(i-1)}$, then ${\rm gr\,}G$ has a natural Lie conformal algebra structure, called the {\it associated graded conformal algebra} of $G$.
\end{defi}
\begin{rema}\label{rema-gc-N}\rm\begin{itemize}\parskip-3pt
\item[(1)]
The ${\rm gr\,}gc_1$ is the Lie conformal algebra with free $\C[\p]$-generating set  $\{J_i\,|\,-1\le i\in\Z\}$, and the $\l$-bracket \vs{-5pt}(cf.~\eqref{gc1-bracket-3})
\begin{equation}\label{BL-con}
[J_{i_{\l}} J_j]=\big((i+j+2)\l+(i+1)\p\big) J_{i+j}\mbox{ \ for \ }i,j\ge-1
.\end{equation}
It is straightforward to verify that ${\rm gr\,}gc_1$ is a simple conformal algebra, whose
corresponding formal distribution Lie algebra is  a well-known Block
type Lie algebra studied in [\ref{S2}], thus we refer this conformal
algebra to as a {\it Block type conformal algebra}.
\item[(2)]
Note from \eqref{gc-N}  that when $N>1$, the filtration of $gc_N$ has to be taken as
\begin{eqnarray}\label{gc-N-fil}\!\!
\{0\}\!=\!gc_N^{(-1)}\!{\sc\!}\subset\! gc_N^{(0)}\!{\sc\!}\subset\!\cdots\!\subset\! gc_N
\mbox{ with }
gc_N^{(n)}\!=\!\mbox{span}\{J^i_A\,|\,A\!\in\! gl_N,\,0\!\le\! i\!\leq\! n\}
\mbox{ for } n\!\geq\!0,
\end{eqnarray}
such \vs{-5pt}that for $i,j,\in\Z_+,\,A,B\in gl_N$,
\begin{equation}\label{gcN-bracket-3}
[J^i_{A_{\l}} J^j_B]\equiv J^{i+j}_{[A,B]}
\mbox{~~(~mod~}gc_N^{(i+j-1)}), \mbox{ where }[A,B]=AB-BA.
\end{equation}
Thus the filtration of $gc_N$ for $N\!>\!1$ is different from that of $gc_1$. This is why we have to treat $gc_1$ separately (if we use the filtration \eqref{gc-N-fil} for $gc_1$, then we obtain from \eqref{gcN-bracket-3}  that the associated graded conformal algebra of this filtration is simply nothing but trivial).
\end{itemize}\end{rema}

The main problem to be addressed in this paper is to
determine filtered Lie conformal algebras whose associated graded
conformal algebras are isomorphic to that of general conformal
algebra $gc_1$.
\vskip7pt

\cl{\bf\S3. \  Some second cohomology groups
}\setcounter{section}{3}\setcounter{equation}{0}\setcounter{theo}{0} \vs{5pt}
To prove Theorem \ref{theo}(1), we need some results on the second cohomology groups. We refer to \cite{BKV} for definition of conformal cohomology.
\def\fg{\Vir{\sc\!}\ltimes{\!} M_{a,0}}Let $\Vir=\C[\ptl]L$  and $M_{\Delta,\alpha}=\C[\ptl]v$
be respectively the Virasoro conformal algebra and its module defined in
\eqref{Vir} and \eqref{M-a-b}.
By \cite[Theorem 7.2(3)]{BKV}, we  \vs{-7pt}have
\begin{eqnarray}\label{H-2-}
&\!\!\!\!\!\!&{\rm dim\,}H^2(\Vir,M_{\Delta,\alpha})=\left\{\begin{array}{ll}
2&\mbox{if \ }\Delta=-1,0\mbox{ and }\alpha=0,
\\[4pt]
1&\mbox{if \ }\Delta=-6,-4,1\mbox{ and }\alpha=0,
\\[4pt]0&\mbox{otherwise}.\end{array}\right.
\\
\label{H-2-1}
&\!\!\!\!\!\!&H^2(\Vir,M_{1,0})=\C\phi,\mbox{ where }\phi_{\l_1,\l_2}(L,L)=\l_1-\l_2.
\end{eqnarray}

Let $\fg=\C[\ptl]L\oplus\C[\ptl]J$ be the semi-direct sum of the Virasoro conformal \mbox{algebra} $\Vir=\C[\ptl]L$ and
its module $M_{a,0}=\C[\ptl]J$. Then a  $\Vir$-module $M_{b,0}=\C[\ptl]v$
\mbox{becomes} a \mbox{$\fg$-module} with the trivial $\l$-action of $J$.
As in \cite[\S7.2]{BKV},
for any $2$-cochain \mbox{$\psi\!\in\! C^2(\fg,M_{b,0})$,} there are $3$ unique
polynomials $P_{JJ}(\l_1,\l_2)$, $P_{\l_1,\l_2},$ $P_{LL}(\l_1,\l_2)$  such that $P_{JJ}(\l_1,\l_2),$ $P_{LL}(\l_1,\l_2)$ are skew-symmetric, and
\begin{equation}\label{CoCy}
\psi_{\l_1,\l_2}(J,J)\!\equiv\! P_{JJ}(\l_1,\l_2)v,\ \ \
\psi_{\l_1,\l_2}(L,J)\!\equiv\! P_{\l_1,\l_2}v,\ \ \
\psi_{\l_1,\l_2}(L,L)\!\equiv\! P_{LL}(\l_1,\l_2)v,\!\!\!\end{equation}
where ``\,$\equiv$\,'' means ``\,equality under modulo $\ptl+\l_1+\l_2$\,''.
\begin{lemm}\label{co-lemm1}
Suppose $\psi\in C^2(\fg,M_{b,0})$ is  a $2$-cocycle. Then for some $c_0\in\C$,\begin{equation}\label{co-JJ}
P_{JJ}(\l_1,\l_2)=\Big\{\begin{array}{ll}c_0(\l_1-\l_2)&\mbox{if \ }b=2a-2,\\[4pt]
 0&\mbox{otherwise}.\end{array}\end{equation}
 \end{lemm}
\noindent{\it Proof.~}~From the definition of differential operator $d$, and using $[J_\l J]=0$ and $J_\l v=0$, we obtain
\begin{eqnarray}\label{A==0}
0&\!\!\!=\!\!\!&(d\psi)_{\l_1,\l_2,\l_3}(L,J,J)
=
L_{\l_1}\psi_{\l_2,\l_3}(J,J)-\psi_{\l_1+\l_2,\l_3}([L_{\l_1}J],J)+
\psi_{\l_1+\l_3,\l_2}([L_{\l_1}J],J)\nonumber\\&\!\!\!\equiv\!\!\!&
(-\l_1-\l_2-\l_3+b\l_1)P_{JJ}(\l_2,\l_3)v
-(-\l_1-\l_2+a\l_1)P_{JJ}(\l_1+\l_2,\l_3)v\nonumber\\&\!\!\!\!\!\!&
+(-\l_1-\l_3+a\l_1)P_{JJ}(\l_1+\l_3,\l_2)v.
\end{eqnarray}

First assume $a\ne1$.
Letting $\l_3=0$, we obtain
\begin{equation}\label{PJJ=01}
P_{JJ}(\l_1,\l_2)=\frac1{(a-1)\l_1}\Big(((1-b)\l_1+\l_2)P_{JJ}(\l_2,0)
+((a-1)\l_1-\l_2)P_{JJ}(\l_1+\l_2,0)\Big).
\end{equation}
If $a=b=2$, then $P_{JJ}(\l_1,\l_2)=\frac{\l_1-\l_2}{\l_1}(P_{JJ}(\l_1+\l_2,0)-P_{JJ}(\l_2,0))$, using this in \eqref{A==0}, we easily obtain that $P_{JJ}(\l,0)=c_0\l$ for some $c_0\in\C$, and the first case of \eqref{co-JJ} holds. If $a=2\ne b$,
using $P_{JJ}(\l,\l)=0$ in \eqref{PJJ=01}, we immediately obtain that $P_{JJ}(\l,0)=0$, thus
 $P_{JJ}(\l_1,\l_2)=0$ by \eqref{PJJ=01}. Hence, we suppose $a\ne2$. Using $P_{JJ}(\l,\l)=0$, we have $
P_{JJ}(2\l,0)=\frac{b-2}{a-2}P_{JJ}(\l,0)$, which implies $P_{JJ}(\l,0)=P_{JJ}(1,0)\l^m$ is a homogenous polynomial of degree, say $m$, such that $\frac{b-2}{a-2}=2^m$.
Then by \eqref{A==0} and \eqref{PJJ=01}, we obtain that $m=1$,  $P_{JJ}(\l,0)=P_{JJ}(1,0)\l $, and  $P_{JJ}(\l_1,\l_2)=\frac{(a-1)\l_1+(a-1-b)\l_2}{a-1}P_{JJ}(1,0)$. This gives  \eqref{co-JJ}.

Finally assume $a=1$. Letting $\l_3=1$ in \eqref{A==0} gives
\begin{equation}\label{a====1}
P_{JJ}(\l_1,\l_2)=((b-1) \l_1  - \l_2-b-1)R({\l_2})
 + \l_2 R({\l_1 + \l_2-1})
,\end{equation}
where $R(\l)=P_{JJ}(\l,1)$. Using $P_{JJ}(\l,0)+P_{JJ}(0,\l)=0$, and using  \eqref{a====1} in \eqref{A==0} with  $\l_2=0$ and $\l_1=\l_3=\l$, we obtain respectively
\begin{eqnarray}\label{pppp0}
&&(\l+b)R(\l)-\l R({\l-1})+((1-b)\l+b)R(0)=0,
\\
&&((b-2)\l-b)R({\l})+\l R({2\l-1})=0.
\label{pppp1}
\end{eqnarray}
If $b\ne0$, then \eqref{pppp1} shows that $R(0)=0$, and from \eqref{pppp0} and \eqref{pppp1}, we can obtain
$R(\l)=0$ for infinite many $\l$'s, thus $R(\l)=0$. If $b=0$, then \eqref{pppp0} gives
$R(1)=0$ and $R(\l)=R(0)(1-\l)$, and we have \eqref{co-JJ} by \eqref{a====1}.
\hfill$\Box$\vs{5pt}
\def\d{\delta}

Now we determine the polynomial $P_{\l_1,\l_2}$ below. Similar to \eqref{A==0}, we have
\begin{eqnarray}\label{JA==0}\!\!\!\!
0&\!\!\!=\!\!\!&(d\psi)_{\l_1,\l_2,\l_3}(L,L,J)
\nonumber\\&\!\!\!=\!\!\!&
L_{\l_1}\psi_{\l_2,\l_3}(L,J)-L_{\l_2}\psi_{\l_1,\l_3}(L,J)\nonumber\\&\!\!\!\!\!\!&
-\psi_{\l_1+\l_2,\l_3}([L_{\l_1}L],J)+
\psi_{\l_1+\l_3,\l_2}([L_{\l_1}J],L)
-\psi_{\l_2+\l_3,\l_1}([L_{\l_2}J],L)
\nonumber\\&\!\!\!\equiv\!\!\!&
((b-1)\l_1-\l_2-\l_3)P_{\l_2,\l_3}v
-(-\l_1+(b-1)\l_2-\l_3)P_{\l_1,\l_3}v\nonumber\\&\!\!\!\!\!\!&
-(\l_1\!-\!\l_2)P_{\l_1+\l_2,\l_3}v\!-\!((a\!-\!1)\l_1\!-\!\l_3)
P_{\l_2,\l_1+\l_3}v\!+\!((a\!-\!1)\l_2\!-\!\l_3)P_{\l_1,\l_2+\l_3}v.
\end{eqnarray}
Denote by $P^{(m)}_{\l_1,\l_2}$ the homogenous part of $P_{\l_1,\l_2}$ of degree $m$.
Then \eqref{JA==0} is satisfied by $P^{(m)}_{\l_1,\l_2}$.
Note that for any polynomial $Q(\l)$, we can define a $1$-cochain $f\in C^1(\fg,M_{b,0})$ by
$f_\l(L)=0,f_\l(J)=Q(\l)v.$
Let $\bar\psi=df$ be the corresponding $2$-coboundary, and set $\psi'=\psi-df$.
Then
\begin{eqnarray}\label{P-df}\!\!\!\!\!\!\!\!&\!\!\!\!\!\!\!\!\!\!\!\!&
\psi'_{\l_1,\l_2}(L,J)=P'_{\l_1,\l_2}v,\mbox{ \ where}\nonumber\\ \!\!\!\!\!\!\!\!&\!\!\!\!\!\!\!\!& P'_{\l_1,\l_2}=P_{\l_1,\l_2}-(-\l_1-\l_2+b\l_1)Q(\l_2)+(-\l_1-\l_2+a\l_1)Q(\l_1+\l_2).
\end{eqnarray}
If we replace $\psi$ by $\psi'$, then $P_{\l_1,\l_2}$ is replaced by $P'_{\l_1,\l_2}$.
Thus, by some suitable choice of $Q(\l)$, we can always suppose
\begin{eqnarray}\label{Assume-P}\begin{array}{ll}
P^{(m)}_{\l,0}=0\mbox{ \ if \ }m\ge1,\,a\ne1, \mbox{ \ or \ }(m,a,b)=(1,1,0),\\[4pt]
P^{(m)}_{\l,-\l}=0\mbox{ \ if \ }m\ge1,\,a=1,\,b\ne0,\\[4pt]
P^{(m)}_{\l,\l}=0\mbox{ \ if \ }m\ge3,\,a=1,\,b=0.
\end{array}\end{eqnarray}

Assume  $P^{(m)}_{\l_1,\l_2}\ne0$. For $m=0,1,2$, one can directly check that
for some $c_i\in\C$,
\begin{eqnarray}\label{P-m=?}
P^{(m)}_{\l_1,\l_2}=\left\{
\begin{array}{ll}
c_1&\mbox{if \ }m=0,\,a=b,\\[4pt]
c_2\l_1&\mbox{if \ }m=1,\,a=b,\\[4pt]
c_3\l_2&\mbox{if \ }(m,a,b)=(1,1,0),\\[4pt]
c_4\l_1\l_2&\mbox{if \ }m=2,\,a=1,\,b\ne0,\\[4pt]
c_5\l_1^2+c_6\l_1\l_2&\mbox{if \ }(m,a,b)=(2,1,0).
\end{array}\right.\end{eqnarray}

From now on, we suppose $m\ge3$.
First assume $a\ne 1,b$. Putting $\l_2=\l_3=0$ in \eqref{JA==0}, we obtain $(a-1)P_{0,\l_1}=(b-1)P_{0,0}$, thus $P_{0,\l_2}=0$. This together with \eqref{Assume-P} proves $P^{(m)}_{\l_1,\l_2}$ is divided by $\l_1\l_2$,
 and we can suppose $P^{(m)}_{\l_1,\l_2}=\l_1\l_2P'_{\l_1,\l_2}$ for some polynomial $P'_{\l_1,\l_2}$.
Now putting $\l_3=0$ in \eqref{JA==0} gives $\l_2P'_{\l_1,\l_2}=\l_1P'_{\l_2,\l_1}$. Thus $P'_{\l_1,\l_2}$ is divided by $\l_1$. Therefore, we can \vs{-5pt}suppose \begin{equation}\label{S-sym}
P^{(m)}_{\l_1,\l_2}=\l_1^2\l_2 S_{\l_1,\l_2} \vs{-5pt},\end{equation}
 where $S_{\l_1,\l_2}$ is a homogenous symmetric polynomial of degree $m-3$.
Write $S_{\l_1,\l_2}$ as
$S_{\l_1,\l_2}=\sum_{i=0}^{m-3}s_{i}\l_1^{m-3-i}\l_2^i$. Comparing the coefficients of $\l_1\l_2^{m-1-j}\l_3^{j+1}$ on the both
sides of \eqref{JA==0}, we have $(m-1+b-a)s_{j}=0$ for $0\leq j\leq m-3$.
Thus $S_{\l_1,\l_2}=0$ if $m\ne a-b+1$.
Now suppose $m=a-b+1$.
For $m=3,4,5,6,7$, one can easily  \vs{-5pt}solve
\begin{equation}\label{P-S---}\!
P^{(m)}_{\l_1,\l_2}\!=\!
\left\{\!\!\begin{array}{ll}c_7\l_1^2\l_2
&\!\!\!\mbox{if  }m\!=\!a\!-\!b\!+\!1\!=\!3,\\[4pt]
c_8\l_1^2\l_2(\l_1+\l_2)&\!\!\!\mbox{if  }m\!=\!a\!-\!b\!+\!1\!=\!4,\\[4pt]
\l_1^2\l_2(c_9(\l_1^2+\l_2^2)+c_{10}\l_1\l_2)&\!\!\!\mbox{if  }m\!=\!a\!-\!b\!+\!1\!=\!5,\\[4pt]
c_{11}\l_1^2\l_2(\l_1^3\!+\!\l_2^3\!+\!{\frac{20}{7}}\l_1\l_2(\l_1\!+\!\l_2))&\!\!\!\mbox{if  }(m,a,b)\!=\!(6,5,0),\\[4pt]
\l_1^2\l_2\big(c_{12} (\l_1^5\! +\!\l_2^5)\! +\! c_{13} (\l_1^4 \l_2 \!+\! \l_1\l_2^4)\! +\!
c_{14}(\l_1^3\l_2^2\! +\! \l_1^2\l_2^3)\big)&\!\!\!\mbox{if }m\!=\!a\!-\!b\!+\!1\!=\!7 \vs{-5pt},
\end{array}\right.\!\!\!\!\!\end{equation}
for some $c_i\in\C$, with the following  \vs{-5pt}conditions,
\begin{eqnarray}\label{P-S---con}\begin{array}{ll}
(10  + 3 b) c_9=( 5+ 2 b) c_{10}\mbox{ if }m\!=\!5,\\[4pt]
(9\!+\!2b)c_{13}\!=\!(33\!+\!5b)c_{12},\ (9\!+\!2b)c_{14}\!=\!(51\!+\!16b)c_{12},\ b\!=\!\frac12(-5\pm\sqrt{19})\mbox{ if }m\!=\!7 \vs{-5pt}.\end{array}
\end{eqnarray}

Now suppose $m=a-b+1\ge8$.
By comparing the coefficients of $\l_1^2\l_2^{a-b-j}\l_3^{j}$ on the both sides of \eqref{JA==0}, we  \vs{-7pt}obtain
\begin{eqnarray}\label{S-1}
(a\!-\!b\!-\!2\!-\!j)(a\!-\!b\!+\!1\!-\!j)s_{j-1}+(j\!+\!1)(j\!-\!2a\!+\!2)s_{j}
=2\big(\big(^{a-b}_{\,\,j}\big)-a\big(^{a-b-1}_{\,\ \ j}\big)\big)s_0,\\
\label{S-2}
(a-b-j)(-a-b+1-j)s_{j-1}+(j-1)(j+2)s_{j}=2\big(\big(^{a-b}_{\,j+1}\big)
-a\big(^{a-b-1}_{\,\ \ j}\big)\big)s_0 \vs{-7pt},
\end{eqnarray} for $1\leq j\leq a-b-2$, where \eqref{S-2} is obtained from \eqref{S-1} by symmetry of $S_{\l_1,\l_2}$. Similarly, by comparing the coefficients of $\l_1^3\l_2^{a-b-j}\l_3^{j-1}$, we  \vs{-5pt}have
\begin{eqnarray}\label{S-3}\!\!\!\!\!\!\!\!\!\!\!\!\!\!&\!\!\!\!\!\!\!\!\!\!\!\!\!\!&\!\!\!\!\!\!\!
(a\!-\!b\!-\!3\!-\!j)(a\!-\!b\!+\!1\!-\!j)
(a\!-\!b\!+\!2\!-\!j)s_{j-2}\!+\!j(j\!+\!1)(j\!-\!3a\!+\!2)
s_{j}
\nonumber\\
\!\!\!\!\!\!\!\!\!\!\!\!\!\!&\!\!\!\!\!\!\!\!\!\!\!\!\!\!&\!\!\!\!\!\!\!
\phantom{(a\!-\!b\!-\!3\!-\!j)(a\!-\!b\!+\!1\!-\!j)
(a\!-\!b\!+\!2\!-\!j)s_{j-2}\!+\!j(j\!+\!1)}
\!=\!6\big(\big(^{a-b-1}_{\ j-1}\big)\!-\!a\big(^{a-b-2}_{\ j-1}\big)\big)s_1,
\\
\label{S-4}
\!\!\!\!\!\!\!\!\!\!\!\!\!\!&\!\!\!\!\!\!\!\!\!\!\!\!\!\!&\!\!\!\!\!\!\!
(j\!+\!b\!-\!a)(a\!-\!b\!+\!1\!-\!j)
(2a\!+\!b\!+\!j\!-\!2)s_{j-2}\!+\!(j\!-\!3)(j\!+\!1)
(j\!+\!2)s_{j}\nonumber\\
\!\!\!\!\!\!\!\!\!\!\!\!\!\!&\!\!\!\!\!\!\!\!\!\!\!\!\!\!&\!\!\!\!\!\!\!
\phantom{(a\!-\!b\!-\!3\!-\!j)(a\!-\!b\!+\!1\!-\!j)
(a\!-\!b\!+\!2\!-\!j)s_{j-2}\!+\!j(j\!+\!1)}
\!=\!6\big(\big(^{a-b-1}_{\,\ j}\big)\!-\!a\big(^{a-b-2}_{\,j-1}\big)\big)s_1,
\end{eqnarray} for $3\leq j\leq a-b-2$. Note that we can recursively use \eqref{S-2} to solve $s_j$ in terms of $s_0$ and $s_1$ for $j\ge2$.
If $a=\frac32$, then from \eqref{S-1}, we obtain $s_0=0$, and one can then easily check from \eqref{S-1}--\eqref{S-4} that $s_j=0$ for all $j$.
Thus suppose $a\ne\frac32$. Use \eqref{S-1} to solve $s_1$ in term of $s_0$, we can then obtain $s_j$ in term of $s_0$. Taking $j=2$ in \eqref{S-1} and $j=4$ in \eqref{S-4}, one immediately obtain $s_0=0$. Thus $S_{\l_1,\l_2}=0$.

Now assume $a=b$.  By \eqref{Assume-P}, we can suppose
$P^{(m)}_{\l_1,\l_2}=\l_2 R_{\l_1,\l_2}$ if $a\ne1$ or
$P^{(m)}_{\l_1,\l_2}=(\l_1+\l_2) R_{\l_1,\l_2}$ if $a=1$, for some homogenous polynomial $R_{\l_1,\l_2}$ of degree $m-1$.
Then using \eqref{JA==0} and discussing as above, we can prove $P^{(m)}_{\l_1,\l_2}=0$.
%
%
%

Next assume $a=1,\,b=0$. As above, we can suppose $P^{(m)}_{\l_1,\l_2}=(\l_1-\l_2)R_{\l_1,\l_2}$  for some homogenous polynomial $R_{\l_1,\l_2}$ of degree $m-1$,
from which we can deduce that  $P^{(m)}_{\l_1,\l_2}=0$.
%
%
%
%

Finally assume $a=1,\, b\ne0,1$.
 As before, we can suppose 
 $P^{(m)}_{\l_1,\l_2}=(\l_1+\l_2)R_{\l,\l_2}$ for some homogenous polynomial $R_{\l_1,\l_2}$ of degree $m-1$.
For $m=3,4,5,6$, one can easily solve
 \begin{equation}\label{P-S---a=1}
P^{(m)}_{\l_1,\l_2}\!=\!
\left\{\!\!\begin{array}{ll}c_{15}\l_1^2(\l_1+\l_2)
&\mbox{if  }(m,a,b)\!=\!(3,1,-1),\\[4pt]
c_{16}\l_1^2\l_2(\l_1+\l_2)&\mbox{if  }
(m,a,b)\!=\!(4,1,-2),\\[4pt]
\l_1^2\l_2^2(\l_1+\l_2)&\mbox{if  }
(m,a,b)\!=\!(5,1,-3),\\[4pt]
c_{17}\l_1^2\l_2(\l_1+\l_2)(\l_1^2+\l_1\l_2+7\l_2^2)
&\mbox{if  }(m,a,b)\!=\!(6,1,-4).
\end{array}\right.\end{equation}
If $m\ge7$, then similar to the arguments after \eqref{P-S---con}, we obtain $P^{(m)}_{\l_1,\l_2}=0$.

Note that every $2$-cochain $\psi\in C^2(\fg,M_{b,0})$ can be restricted to a $2$-cochain $\psi':=\psi|_{\Vir\times\Vir}\in C^2(\Vir,M_{b,0})$,
and conversely, every $2$-cochain $\psi'\in C^2(\Vir,M_{b,0})$ can be extended to a $2$-cochain $\psi\in C^2(\fg,M_{b,0})$ by taking the corresponding polynomials $P_{JJ}(\l_1,\l_2)$, $P_{\l_1,\l_2}$ to be zero. From this, we obtain an embedding
$H^2(\Vir,M_{b,0})\!\to\! H^2(\fg,M_{b,0}).$ Thus, we can regard $H^2(\Vir,M_{b,0})$ as a subspace of $H^2(\fg,M_{b,0})$. Now we can state the main result in this section.

\begin{theo}\label{Prop1} We have
\begin{eqnarray}\label{PROP}
\!\!\!\!\!\!\!\!&\!\!\!\!\!\!\!\!&{\rm dim\,}H^2(\fg,M_{b,0})={\rm dim\,} H^2(\Vir,M_{b,0})+\d_{b,2a-2}+
\tau_{a,b},\mbox{ where }\\
\nonumber\!\!\!\!\!\!\!\!&\!\!\!\!\!\!\!\!&
\tau_{a,b}=\left\{\begin{array}{ll}
3&\mbox{if \ }(a,b)\!=\!(1,0),\\[4pt]
2&\mbox{if \ }a\!=\!b,\mbox{ or \ }a\!=\!1,\ b\!=\!-3,-4,-5,-6,\\[4pt]
1&\mbox{if \ }a\!=\!1,\ b\!\ne\!1,0,-3,-4,-5,-6,\mbox{ or \ }a\!\ne\!1,\ a\!-\!b\!=\!2,3,4,
\\[2pt]&\mbox{or \ }a\!=\!5,\ b\!=\!0,
\mbox{ or \ }a\!=\!6\!+\!b,\ b\!=\!\frac12(-5\!\pm\!\sqrt{19}),\\[4pt]
0&\mbox{otherwise.}
\end{array}\right.
\end{eqnarray}
Furthermore, every $2$-cocycle $\psi'\in C^2(\fg,M_{b,0})$ is equivalent to a $2$-cocycle $\psi$ such that
the corresponding polynomial $P_{JJ}(\l_1,\l_2)$ defined in \eqref{CoCy} has the form in \eqref{co-lemm1}, and the homogenous part $P^{(m)}_{\l_1,\l_2}$ of  $\,P_{\l_1,\l_2}$ has the form in    \eqref{P-m=?}, \eqref{P-S---} or \eqref{P-S---a=1}.
\end{theo}
\vs{5pt}

\cl{\bf\S4. \  Proof of Theorem \ref{theo}(1)
}\setcounter{section}{4}\setcounter{equation}{0}\setcounter{theo}{0} \vs{5pt}

Now in order to prove the main theorem, we need some preparations.
Since $G$ is a filtered Lie conformal algebra and ${\rm gr\,}G\cong
{\rm gr\,}gc_1 $, by $(\ref {gc1-bracket-3})$, we can suppose
\begin{eqnarray}\label{J_-1_-1}
\!\!\!\!\!\!\!\!\!\!\!\!&\!\!\!\!\!\!\!\!\!\!\!\!\!\!\!&
[J{_{-1}}_{\l}J_{-1}]=0\ ,
\\\!\!\!\!\!\!\!\!\!\!\!\!&\!\!\!\!\!\!\!\!\!\!\!\!\!\!\!&
[J{_{-1}}_{\l}J_0]=\l J_{-1}
\ ,\ \label{J_-1_0}
\\\!\!\!\!\!\!\!\!\!\!\!\!&\!\!\!\!\!\!\!\!\!\!\!\!\!\!\!&
[J{_0}_\l J_0]=(2\l+\p)J_0+g_1(\l,\p)J_{-1} \ ,\label{J_0_0}
\\\!\!\!\!\!\!\!\!\!\!\!\!&\!\!\!\!\!\!\!\!\!\!\!\!\!\!\!&
[J{_{-1}}_{\l}J_1]=2\l J_{0}+g_2(\l,\p)J_{-1}\ ,\label{J_-1_1}
\\\!\!\!\!\!\!\!\!\!\!\!\!&\!\!\!\!\!\!\!\!\!\!\!\!\!\!\!&
[J{_0}_{\l}J_1]=(3\l+\p)J_1+g_3(\l,\p) J_{0}+g_4(\l,\p)J_{-1}\ ,\label{J_0_1}
\\\!\!\!\!\!\!\!\!\!\!\!\!&\!\!\!\!\!\!\!\!\!\!\!\!\!\!\!&
[J{_{-1}}_{\l}J_2]=3\l J_1+h_1(\l,\p) J_{0}+h_2(\l,\p)J_{-1}\ ,\label{J_-1_2}
\\\!\!\!\!\!\!\!\!\!\!\!\!&\!\!\!\!\!\!\!\!\!\!\!\!\!\!\!&
[J{_{0}}_{\l}J_2]=(4\l+\p) J_2+h_3(\l,\p) J_{1}+h_4(\l,\p) J_{0}+h_5(\l,\p)J_{-1}\ ,\label{J_0_2}
\\\!\!\!\!\!\!\!\!\!\!\!\!&\!\!\!\!\!\!\!\!\!\!\!\!\!\!\!&
[J{_{1}}_{\l}J_1]=2(2\l+\p) J_2+f_1(\l,\p) J_{1}+f_2(\l,\p) J_{0}+f_3(\l,\p)J_{-1}\ ,\label{J_1_1}
\end{eqnarray}
where $g_i(\l,\p)$ for $1\leq i\leq 4$,  $h_i(\l,\p)$ for $1\leq
i\leq 5$ and $f_i(\l,\p)$ for $1\leq i\leq 3$  are all polynomials
of $\l$ and $\p$. Then our aim is to determine all these polynomials
of $\l$ and $\p$ by making use of Theorem \ref{Prop1} and by computations.
\begin{lemm}\label{lemm1} In $(\ref{J_0_0})$, \RED{by re-choosing the generator $J_0$, we
can suppose}
$ 
g_1(\l,\p)=0.
$ 
\end{lemm}
\noindent{\it Proof.~}~Noting from  \eqref{J_-1_0}, we see $P_{LL}(\l_1,\l_2):=g_1(\l_1,-\l_1-\l_2)$ is a skew-symmetric
 polynomial, thus we can define a $2$-cochain $\psi\in C^2(\Vir,M_{1,0})$ with $\psi_{\l_1,\l_2}(L,L)=g_1(\l_1,-\l_1-\l_2)$, which is in fact a $2$-cocycle by Jacobi identity.
Note that for any polynomial $p(\ptl)$, if we replace the generator $J_0$ by $J'_0=J_0-p(\ptl)J_{-1}$, it is equivalent to replacing the $2$-cocycle $\psi$ by $\psi'=\psi-d\phi$, where $\phi\in C^1(\Vir,M_{1,0})$ is the $1$-cochain defined by $\phi_\l(L)=p(-\l)$. Thus by \eqref{H-2-1}, we can suppose $g_1(\l_1,-\l_1-\l_2)=a_0(\l_1-\l_2)$ for some $a_0\in\C$, i.e.,
$ 
g_1(\l,\p)=a_0(\ptl+2\l).
$ 
Applying the operator ${J_0}_{\mu}$ to \eqref{J_-1_1} and \eqref{J_0_1}, using the Jacobi identity and comparing the coefficients of $J_{0}$ and $J_{-1}$ respectively, we obtain
\begin{eqnarray*}\!\!\!\!\!\!\!\!\!\!\!\!&\!\!\!\!\!\!\!\!\!\!\!\!\!\!\!&
\ \ \ \ \ \ \ \ (3\m+\l+\p)g_3(\l,\p)+(2\l+\p)g_3(\m,\l+\p)-(\l-\m)g_3(\l+\m,\p)
\nonumber\\\!\!\!\!\!\!\!\!\!\!\!\!&\!\!\!\!\!\!\!\!\!\!\!\!\!\!\!&
\ \ \ \ \ \ \ \ =(3\l+\m+\p)g_3(\m,\p)+(2\m+\p)g_3(\l,\p+\m)+2a_0(\l+\m)(\l-\m).
\\[4pt]
\!\!\!\!\!\!\!\!\!\!\!\!&\!\!\!\!\!\!\!\!\!\!\!\!\!\!\!&\ \ \ \ \ \ \ \
2\mu{\ssc\,} a_0(\p+2\l)+(\l+\p)g_2(\mu,\l+\p)+\mu g_2(\l+\mu,\p)
\nonumber\\\!\!\!\!\!\!\!\!\!\!\!\!&\!\!\!\!\!\!\!\!\!\!\!\!\!\!\!&\ \ \ \ \ \ \ \
=(3\l+\mu+\p)g_2(\mu,\p)+\mu g_3(\l,\mu+\p).
\end{eqnarray*}
To prove $a_0=0$, we only need to suppose $g_3(\l,\p),g_2(\l,\p)$ are homogenous polynomials of degree $1$ (since coefficients of terms of other degrees do not contribute). Then one can immediately check from the above two equations that  $a_0=0$. The lemma follows.
\hfill$\Box$\vspace{5pt}

Similar to the proof of Lemma \ref{lemm1}, from \eqref{J_0_1},
we can use $g_3(\l,\p)$
to define a $2$-cocycle $\psi\in$\linebreak $C^2(\Vir\ltimes M_{3,0},M_{2,0})$ such that the corresponding polynomials defined in \eqref{CoCy} have the forms,
\begin{equation*}
P_{LL}(\l_1,\l_2)=P_{JJ}(\l_1,\l_2)=0, \mbox{ \ and \ }
P_{\l_1,\l_2}=g_3(\l_1,-\l_1-\l_2).
\end{equation*}
Thus, using Theorem \ref{Prop1},  by replacing $J_1$ by $J'_1=J_1+p(\p)J_0$ for some polynomial $p(\p)$, we can suppose $g_3(\l,\p)=0$.
Similarly, $g_4(\l,\p)$ defines a $2$-cocycle in $C^2(\Vir\ltimes M_{3,0},M_{1,0})$, using Theorem \ref{Prop1} and the first case of \eqref{P-S---}, by replacing $J_1$ by $J'_1=J_1+p(\p)J_{-1}$ for some polynomial $p(\p)$ (this replacement does not affect $g_3(\l,\p)$), we can suppose
\begin{equation*}
g_4(\l,\p)=a_1\l^2(\l+\p)\mbox{ \ for some \ $a_1\in\C$.}\end{equation*}

\begin{lemm}\label{lemm1-2} In $(\ref {J_0_1})$, there exists some $a_2\in\C$ such that
$ 
g_2(\l,\p)=a_2\l(\l-\p).
$ 
\end{lemm}
\noindent{\it Proof.~}~Comparing the
coefficients of $J_{-1}$ on the both sides of the Jacobi identity
$[J{_0}_{\l}[J{_{-1}}_\m J_1]]=
[[J{_0}_\l J_{-1}]_{\l+\m} J_1]+[J{_{-1}}_\m[J{_0}_\l J_1]],$  by
$(\ref{J_-1_0})$--$(\ref{J_0_1})$, Lemma \ref{lemm1} and $g_3(\l,\p)=0$, we
have
\begin{equation}\label{g2-1}
(\l+\p)g_2(\m,\l+\p)-(3\l+\m+\p)g_2(\m,\p)+\m g_2(\l+\m,\p)=0.
\end{equation}
Taking $\p=0$ gives
\begin{equation}\label{g2-1+1-1}
g_2(\mu,\l)=\frac1{\l}
\big((3\l+\mu)g_2(\mu,0)-\mu g_2(\l+\mu,0)\big).
\end{equation}
Using this in \eqref{g2-1} and taking $\p=\mu=-\l$, we obtain
$g_2(2\l,0)=4g_2(\l,0)+3g_2(0,0)$. Inductively, we obtain that
$ 
g_2(x\l,0)=x^2g_2(\l,0)+(x^2-1)g_2(0,0)
$ 
holds for all $x=2^k,\,k=1,2,...$, thus holds for all $x\in\C$ since $g_2(\l,\mu)$ is a polynomial.
Taking $x=0$, we obtain $g_2(0,0)=0$. Thus $g_2(x,0)=a_2x^2$, where $a_2=g_2(1,0)$.
Using this in \eqref{g2-1+1-1}, we obtain $g_2(\l,\p)=a_2\l(\l-\p)$. \hfill$\Box$\vs{5pt}

For later convenience (in order to let other polynomials have some suitable forms),
we respectively re-denote $a_2$ by $\frac{c}2$, $a_1$ by $-\frac{b}2$ and replace $J_1$ by $J'_1=J_1+\frac{c}2\p J_0+\frac{b}{2}\p^2J_{-1}$, so that $g_2(\l,\p)$, $g_3(\l,\p)$ and $g_4(\l,\p)$ have the following forms,
\begin{equation}\label{g-2-3-4}
g_2(\l,\p)=c\l^2, \ \ \ g_3(\l,\p)=c\l^2, \ \ \ g_4(\l,\p)=b\l^2\p
\mbox{ \ \ for some \ $b,c \in\C$.}\end{equation}
Similarly, by \eqref{J_0_2}, $h_3(\l,\p),\,h_4(\l,\p)$ and $h_5(\l,\p)$ define {\def\fg{{\mathfrak g}}$2$-cocycles $\psi_3\in C^2(\fg,M_{3,0})$,
$\psi_4\in C^2(\fg,M_{2,0})$ and $\psi_5\in C^2(\fg,M_{1,0})$ respectively with $\fg=\Vir\ltimes M_{4,0}$.}
Thus using Theorem \ref{Prop1}, by replacing $J_2$ by $J_2'=J_2+p_0(\p)J_1+p_1(\p)J_0+p_2(\p)J_{-1}$ for some polynomials $p_0(\p)$, $p_1(\p)$, $p_2(\p)$, we can suppose  for some $a_3,t\in\C$,
\begin{equation}\label{h3-h4-h5}
h_3(\l,\p)=0,\ \ \ h_4(\l,\p)=a_3\l^2(\l+\p),\ \ \
h_5(\l,\p)=t\l^2(\l+\p)\p.\end{equation}

\begin{lemm}\label{lemm2} In $(\ref {J_-1_2})$, we have
$ 
h_1(\l,\p)=c\l(\l-2\p).
$ 
\end{lemm}
\noindent{\it Proof.~}~Comparing the coefficients of $J_{-1}$ and $J_0$ on the both sides of the Jacobi identities \linebreak $[J{_{-1}}_{\l}[J{_{-1}}_\m J_2]]=[J{_{-1}}_\m[J{_{-1}}_\l J_2]]$ and
$[J{_{-1}}_{\l}[J{_{0}}_\m J_2]]=[J{_{-1}}_\l J{_{0}}]_{\l+\mu} J_2]]+[J{_{0}}_\m[J{_{-1}}_\l J_2]]$ respectively, we deduce
\begin{eqnarray}
\label{h1-1}\!\!\!\!\!\!\!\!&\!\!\!\!\!\!\!\!&
3\m g_2(\l,\p)+\l h_1(\m,\l+\p)=3\l g_2(\m,\p)+\m h_1(\l,\m+\p),\\
\!\!\!\!\!\!\!\!&\!\!\!\!\!\!\!\!&
(\l+4\m+\p)h_1(\l,\p)+2\l h_3(\m,\l+\p)\nonumber\\
\!\!\!\!\!\!\!\!&\!\!\!\!\!\!\!\!&
=\l h_1(\l+\m,\p)+3\l g_3(\m,\p)+(2\m+\p)h_1(\l,\m+\p)\label{h1-1'}.
\end{eqnarray}
Set $\m=\p=0$ in $(\ref{h1-1})$, by $(\ref{g-2-3-4})$, we
 get
$
\l h_1(0,\l)=2\l g_1(0,0)=0$, i.e., $h_1(0,\l)=0$.
Taking $\p=0$ in $(\ref{h1-1'})$, by $(\ref{g-2-3-4})$ and $(\ref{h3-h4-h5})$, we immediately obtain
\begin{equation}\label{h1-7}
h_{1}(\l,\m)=\frac{1}{2\m}\big((\l+4\m)h_1(\l,0)-\l h_1(\l+\m,0)-3c\l\m^2\big).
\end{equation}
Using this and $(\ref{g-2-3-4})$ in $(\ref{h1-1})$, we obtain
\begin{eqnarray*}
\!\!\!\!\!\!\!\!&\!\!\!\!\!\!\!\!&
\l(\m+\p)(4\l+\m+4\p)h_1(\m,0)+\l\m(\l-\m) h_1(\l+\p+\m,0)\nonumber\\
\!\!\!\!\!\!\!\!&\!\!\!\!\!\!\!\!&
=\m(\l+\p)(\l+4\m+4\p) h_1(\l,0)+3c\l \m (\m-\l)(\l+\p)(\m+\p).
\end{eqnarray*}
Setting $\p=-\l-1,\,\m=1$, and using $h_1(0,0)=0$, we obtain $h_1(\l,0)=\l h_1(1,0)+c\l(\l-1)$. Thus $(\ref{h1-7})$ turns into $h_{1}(\l,\p)=c\l(\l-2\p)+\frac{3}{2}\l(h_1(1,0)-c)$. Taking $\l=1$ and $\p=0$ gives $h_1(1,0)=c$. Hence we have the lemma.\hfill$\Box$ \vskip5pt

For later convenience (in order to let other polynomials have some suitable forms),
we re-denote $a_3$ by $\frac{5k}3$, and replace $J_2$ by $J'_2=J_2+ c\p J_1-\frac{k}{3}\p^2 J_0$, so that by \eqref{h3-h4-h5} and Lemma \ref{lemm2}, $h_1(\l,\p)$, $h_3(\l,\p),$ $h_4(\l,\p)$, $h_5(\l,\p)$ have the following forms
\begin{equation}
h_1(\l,\p)=3c\l^2,\ \ h_3(\l,\p)=3c\l^2,\ \ h_4(\l,\p)=k\l^3,\ \ h_5(\l,\p)=t\l^2(\l+\p)\p.\label{h-1-3-4}
\end{equation}
\begin{lemm}\label{lemm6} In $(\ref {J_-1_2})$, we have
\begin{eqnarray}
\label{h2}\!\!\!\!\!\!\!\!&\!\!\!\!\!\!\!\!&
h_2(\l,\p)=\frac{1}{2}(k-b-c^2)\l\p^2-\frac{3}{2}(k-c^2+b)\l^2\p+\frac{1}{2}(k+b+c^2)\l^3.\\
\label{ffff1}\!\!\!\!\!\!\!\!&\!\!\!\!\!\!\!\!&f_1(\l,\p)=-c(2\l\p+\p^2).
\end{eqnarray}
\end{lemm}
\noindent{\it Proof.~} Comparing the coefficients of $J_{-1}$ on the
both sides of $[J{_{0}}_{\l}[J{_{-1}}_\m J_2]]$
$=[[J{_0}_{\l}J_{-1}]_{\l+\m}J_2]+[J{_{-1}}_\m[J{_{0}}_\l J_2]]$, we
get
\begin{eqnarray*}\!\!\!\!\!\!\!\!\!\!\!\!&\!\!\!\!\!\!\!\!\!\!\!\!\!\!\!&
\ \ \ \ \ \ \ \ 3\m g_4(\l,\p)-g_2(\m,\p)h_3(\l,\m+\p)-\m h_4(\l,\m+\p)
\nonumber\\\!\!\!\!\!\!\!\!\!\!\!\!&\!\!\!\!\!\!\!\!\!\!\!\!\!\!\!&
\ \ \ \ \ \ \ \ =-\m h_2(\l+\m,\p)+(4\l+\m+\p)h_2(\m,\p)-(\l+\p)h_2(\m,\l+\p).
\end{eqnarray*}
Since we already get $g_4(\l,\p)=b\l^2\p$, $g_2(\l,\p)=c\l^2$, $h_3(\l,\p)=3c\l^2$ and $h_4(\l,\p)=k\l^3$,
the above equation turns \vs{-5pt}into
\begin{eqnarray*}
\ \
\m h_2(\l\!+\!\m,\p)-(4\l+\m+\p)h_2(\m,\p)+(\l+\p)h_2(\m,\l\!+\!\p)
=-3b\l^2\m\p+3c^2\l^2\m^2+k\l^3\m\vs{-5pt}.
\end{eqnarray*}
Setting $\m=0$ gives
$(4\l+\p)h_2(0,\p)=(\l+\p)h_2(0,\l+\p)$,
which means $h_2(0,\p)$ is divided by $\l+\p$. Thus $h_2(0,\p)=0$.
Then taking $\l=-\mu$ and $\p=0$ respectively in the above equation, we obtai\vs{-5pt}n
\begin{eqnarray}
\!\!\!\!&\!\!\!\!\!\!\!\!&(3\mu-\p)h_2(\mu,\p)-(\mu-\p)
h_2(\mu,\p-\mu)=-3b\mu^3\p+(3c^2-k)\mu^4,\nonumber\\
\label{SMSMS-2}
\!\!\!\!&\!\!\!\!\!\!\!\!&
h_2(\m,\l)=\frac{1}{\l}
(-\m
h_2(\l+\m,0)+(4\l+\m)h_2(\m,0)
 +3c^2\l^2\m^2+k\l^3\m)\vs{-5pt}.
\end{eqnarray}
Using the second equation in the first equation with $\l=\p=\mu$, we obtain
\begin{eqnarray}\label{h2=====}
h_2(2\mu,0)=5h_2(\mu,0)+\frac{3}{2}(k+b+c^2)\m^3.\end{eqnarray}
If $h_2(\mu,0)$ has degree, say $m$, greater than $3$, then by comparing coefficients of $\mu^m$ in \eqref{h2=====}, we obtain $2^m=5$, a contradiction. Thus $m\le3$, and we can suppose $h_2(\mu,0)=d_3\m^3+d_2\m^2+d_1\mu+d_0$. By \eqref{h2=====}, we immediately get $d_3=\frac{1}{2}(k+b+c^2)$ and $d_2=d_1=d_0=0$. Now using this in \eqref{SMSMS-2}, we obtain \eqref{h2}.

By comparing the coefficients of $J_{0}$ on the both sides of the Jacobi identity \begin{equation}\label{J-111}
[J{_{-1}}_{\l}[J{_{1}}_\m J_1]]=[[J{_{-1}}_{\l}J_{1}]_{\l+\m}J_1]+[J{_{1}}_\m[J{_{-1}}_\l J_1]],\end{equation}
we deduce
\begin{eqnarray*}\!\!\!\!\!\!\!\!\!\!\!\!&\!\!\!\!\!\!\!\!\!\!\!\!\!\!\!&
\ \ \ \ \ \ \ \ (2\m+\l+\p)h_1(\l,\p)+\l g_3(-\m-\p,\p)+\l f_1(\m,\l+\p)
\nonumber\\\!\!\!\!\!\!\!\!\!\!\!\!&\!\!\!\!\!\!\!\!\!\!\!\!\!\!\!&
\ \ \ \ \ \ \ \ =\l g_3(\l+\m,\p)+(\l+\m)g_2(\l,-\l-\m)+(\m+\p)g_2(\l,\m+\p).
\end{eqnarray*}
Taking $\p=0$ and sing  $g_3(\l,\p)=c\l^2$, $g_2(\l,\p)=c\l^2$ and $h_1(\l,\p)=3c\l^2$, we
immediately obtain \eqref{ffff1}.
\hfill$\Box$

\begin{lemm}\label{lemm8} We have
$f_2(\l,\p)=\frac{2b}{5}(2\l\p^2+\p^3),$ and $k=c^2-\frac{7}{5}b.$\end{lemm}
\noindent{\it Proof.~}~Comparing the coefficients of $J_{-1}$ on the both sides of \eqref{J-111}, we deduce
\begin{eqnarray*}\!\!\!\!\!\!\!\!\!\!\!\!&\!\!\!\!\!\!\!\!\!\!\!\!\!\!\!&
\ \ \ \ \ \ \ \ 2(2\m+\l+\p)h_2(\l,\p)+f_1(\m,\l+\p)g_2(\l,\p)+\l f_2(\m,\l+\p)
\nonumber\\\!\!\!\!\!\!\!\!\!\!\!\!&\!\!\!\!\!\!\!\!\!\!\!\!\!\!\!&
\ \ \ \ \ \ \ \ =2\l g_4(\l+\m,\p)-2\l g_4(-\m-\p,\p)-g_2(\l,\m+\p)g_2(-\m-\p,\p)
\nonumber\\\!\!\!\!\!\!\!\!\!\!\!\!&\!\!\!\!\!\!\!\!\!\!\!\!\!\!\!&
\ \ \ \ \ \ \ \ \ \ \ +g_2(\l,-\l-\m)g_2(\l+\m,\p).
\end{eqnarray*}
Using $(\ref{g-2-3-4})$, $(\ref{h-1-3-4})$ and \eqref{ffff1}, we obtain \begin{eqnarray*}\!\!\!\!\!\!\!\!\!\!\!\!&\!\!\!\!\!\!\!\!\!\!\!\!\!\!\!&
\l f_{2}(\m,\l+\p)
=-\l(\l+2\m+\p)\big((k+b-c^2)(\l^2+\p^2)+(-5b+3c^2-3k)\l\p\big).
\end{eqnarray*}
Taking $\p=0$ and $\p=\l$ respectively, we immediately obtain the result.
\hfill$\Box$\vs{5pt}

Now we collect information we obtain in Lemmas $\ref{lemm1}$--$\ref{lemm8}$ as follows \begin{eqnarray}\!\!\!\!\!\!\!\!\!\!\!\!&\!\!\!\!\!\!\!\!\!\!\!\!\!\!\!&
\label{ver-4}\ \ \ \ g_1(\l,\p)=0,\ \ \ \ \ \ \ \ \ g_2(\l,\p)=c\l^2,\ \ \ \ \ \ \ \ \ g_3(\l,\p)=c\l^2,\ \ \ \ \ \ \ \ \ \ \ \ g_4(\l,\p)=b\l^2\p,
\nonumber\\\!\!\!\!\!\!\!\!\!\!\!\!&\!\!\!\!\!\!\!\!\!\!\!\!\!\!\!&
\ \ \ \ h_1(\l,\p)=3c\l^2,\ \ \ \ h_2(\l,\p)=-\frac{6}{5}b\l\p^2+\frac{3}{5}b\l^2\p+(c^2-\frac{1}{5}b)\l^3,\ \ \ \ h_3(\l,\p)=3c\l^2,
\nonumber\\\!\!\!\!\!\!\!\!\!\!\!\!&\!\!\!\!\!\!\!\!\!\!\!\!\!\!\!&
\ \ \ \ h_4(\l,\p)=(-\frac{7}{5}b+c^2)\l^3,\ \ \ \ h_5(\l,\p)=t(\l^2\p^2+\l^3\p),
\nonumber\\\!\!\!\!\!\!\!\!\!\!\!\!&\!\!\!\!\!\!\!\!\!\!\!\!\!\!\!&
\ \ \ \ f_1(\l,\p)=-c(2\l\p+\p^2),\ \ \ \ f_2(\l,\p)=\frac{2}{5}b(\p^3+2\l\p^2).
\end{eqnarray}
Next we should determine $f_3(\l,\p)$. By comparing the coefficients of $J_{-1}$ on the both sides of the Jacobi identity  $ [J{_{0}}_{\l}[J{_{1}}_\m J_1]]=[[J{_{0}}_{\l}J_{1}]_{\l+\m}J_1]+[J{_{1}}_\m[J{_{0}}_\l J_1]]$, we obtain that
\begin{eqnarray}\!\!\!\!\!\!\!\!\!\!\!\!&\!\!\!\!\!\!\!\!\!\!\!\!\!\!\!&
\label{ver-5}\ \ \ \ \ 2(2\m+\l+\p)h_5(\l,\p)+f_1(\m,\l+\p)g_4(\l,\p)+(\l+\p)f_3(\m,\l+\p)
\nonumber\\\!\!\!\!\!\!\!\!\!\!\!\!&\!\!\!\!\!\!\!\!\!\!\!\!\!\!\!&
\ \ \ \ \ =(2\l-\m)f_3(\l+\m,\p)+g_3(\l,-\l-\m)g_4(\l+\m,\p)+g_4(\l,-\l-\m)g_2(\l+\m,\p)
\nonumber\\\!\!\!\!\!\!\!\!\!\!\!\!&\!\!\!\!\!\!\!\!\!\!\!\!\!\!\!&
\ \ \ \ \ \ \ \ +(3\l+\m+\p)f_3(\m,\p)-g_3(\l,\m+\p)g_4(-\m-\p,\p)-g_4(\l,\m+\p)g_2(-\m-\p,\p).
\end{eqnarray}

\begin{lemm}\label{lemm10}We have $f_3(\l,\p)=-bc(\p^4+3\l\p^3+3\l^2\p^2+2\l^3\p)$ and $t=\frac{3}{2}bc$.\end{lemm}
\noindent{\it Proof.~}~Using $(\ref{ver-4})$ in $(\ref{ver-5})$ gives
\begin{eqnarray}\!\!\!\!\!\!\!\!\!\!\!\!&\!\!\!\!\!\!\!\!\!\!\!\!\!\!\!&
\label{f3-1}\ \ \ \ \ (\l+\p)f_3(\m,\l+\p)-(2\l-\m)f_3(\l+\m,\p)-(3\l+\m+\p)f_3(\m,\p)
\nonumber\\\!\!\!\!\!\!\!\!\!\!\!\!&\!\!\!\!\!\!\!\!\!\!\!\!\!\!\!&
\ \ \ \ \ =-\l^2(\l+2\m+\p)\big(bc(\l^2+\l\m+\m^2-3\l\p+\m\p+\p^2)+2t(\p^2+\l\p)\big).
\end{eqnarray}
Taking $\p=0$ gives
\begin{eqnarray}\!\!\!\!\!\!\!\!\!\!\!\!&\!\!\!\!\!\!\!\!\!\!\!\!\!\!\!&
\label{f3-3}\ \ \ \ \ f_3(\m,\l)=\frac{1}{\l}\big((2\l-\m)f_3(\l+\m,0)-(3\l+\m)f_3(\m,0)\big)
\nonumber\\\!\!\!\!\!\!\!\!\!\!\!\!&\!\!\!\!\!\!\!\!\!\!\!\!\!\!\!&
\ \ \ \ \ \ \ \ \ \ \ \ \ \ \ \ \ \ \ -\l(\l+2\m)bc(\l^2+\l\m+\m^2).
\end{eqnarray}
Using this in \eqref{f3-1} with $\p=-\mu=-\l$, we obtain
$
f_3(2\l,0)=-10f_3(\l,0)-9f_3(0,0).$ 
If $f_3(\l,0)$ has degree, say $m$, greater than $0$, then by comparing coefficients of $\l^m$, we obtain $2^m=-10$, a contradiction. Thus $m=0$, and so $f_3(\l,0)=f_3(0,0)=0$. Now using this in \eqref{f3-3}, we obtain that $f_3(\l,\p)=-bc(2\l+\p)\p(\p^2+\l\p+\l^2)$. Using this in \eqref{f3-1} gives $t=\frac{3}{2}bc$. Therefore the lemma holds.\hfill$\Box$\vskip5pt

In order to prove Theorem $\ref{theo}$(1), we need the following lemma whose proof seems to be rather technical.

\begin{lemm}\label{lemm11}In $(\ref{ver-4})$, we have $b=0$.
\end{lemm}

\noindent{\it Proof.~} First we assume $c=0$.
By $(\ref{gc1-bracket-3})$, we can suppose
$$[J{_1}_\l J_2]=(5\l+2\p)J_3+f_4(\l,\p)J_2+f_5(\l,\p)J_1+f_6(\l,\p)J_0+f_7(\l,\p)J_{-1},$$ where, $f_i(\l,\p)$ are polynomials of $\l$ and $\p$ for $4\leq i\leq 7$. From $(\ref{J_a_b_c})$, we can get the Jacobi identity
\begin{equation}\label{b-1}
[J{_1}_{\l}[J{_{1}}_\m J_1]]=[[J{_1}_{\l}J_{1}]_{\l+\m}J_1]+[J{_{1}}_\m[J{_1}_\l J_1]].\end{equation}
Since our purpose is to prove $b=0$,
by comparing the coefficients of $J_{1}$ and $J_{-1}$ on the both sides
of \eqref{b-1} respectively,
we see that the
coefficients of terms of $f_5(\l,\p)$ (resp., $f_7(\l,\p)$)
with degree not equal to $3$ (resp., $5$) do not have relations
 with $b$.
Thus
we can suppose that $f_5(\l,\p)=\sum_{i=0}^{3}a_i\l^{3-i}\p^i$
and $f_7(\l,\p)=\sum_{i=0}^{5}b_i\l^{5-i}\p^i$.
Set $J_3'=J_3+c_1\p^2J_1+c_2\p^4J_{-1}$ for some $c_1$, $c_2\in \C$, then we can choose suitable complex numbers $c_1$ and $c_2$ such that by replacing $J_3$ by $J'_3$, we can suppose $a_3=0$, $b_5=0$. Thus  \begin{equation}\label{b-2} f_5(\l,\p)=\mbox{$\sum\limits_{i=0}^{2}$}a_i\l^{3-i}\p^i,\,\,\,\, f_7(\l,\p)=\mbox{$\sum\limits_{i=0}^{4}$}b_i\l^{5-i}\p^i.\end{equation}
Since $c=0$, by $(\ref{ver-4})$ and Lemma $\ref{lemm10}$, comparing the coefficients of $J_1$ and $J_{-1}$ on the both sides of $(\ref{b-1})$ respectively, we obtain that \begin{eqnarray}\!\!\!\!\!\!\!\!\!\!\!\!&\!\!\!\!\!\!\!\!\!\!\!\!\!\!\!&
\label{b-3}\ \ \ \ \ [J{_1}_\l J_2]=(5\l+2\p)J_3+f_4(\l,\p)J_2+(a_0\l^3+(a_0-b)\l^2\p+a_2\l\p^2)J_1
\nonumber\\\!\!\!\!\!\!\!\!\!\!\!\!&\!\!\!\!\!\!\!\!\!\!\!\!\!\!\!&
\ \ \ \ \ \ \ \ \ \ \ \ \ \ \ \ \ \ +f_6(\l,\p)J_0+(\frac{1}{5}b^2\l^4\p+b_2\l^3\p^2+(b_2+\frac{2}{5}b^2)\l^2\p^3)J_{-1}.
\end{eqnarray}
By $(\ref{J_a_b_c})$, we also have the Jacobi \vs{-5pt}identity
\begin{equation}\label{LALALA}[J{_0}_{\l}[J{_{1}}_\m J_2]]=[[J{_0}_{\l}J_{1}]_{\l+\m}J_2]+[J{_{1}}_\m[J{_0}_\l J_2]]\vs{-5pt}.\end{equation}
 By $(\ref{ver-4})$ (with $c=0$) and $(\ref{b-3})$, we obtain  \begin{eqnarray}\!\!\!\!\!\!\!\!\!\!\!\!&\!\!\!\!\!\!\!\!\!\!\!\!\!\!\!&
\label{b-4}\ \ \ \ \ (5\m+2\l+2\p)[J{_0}_\l J_3]=(5\m+2\l+2\p)(5\l+\p)J_3+l_1(\l,\m,\p)J_2+l_2(\l,\m,\p)J_1
\nonumber\\\!\!\!\!\!\!\!\!\!\!\!\!&\!\!\!\!\!\!\!\!\!\!\!\!\!\!\!&
\ \ \ \ \ \ \ \ \ \ \ \ \ \ \ \ \ \ \ \ \ \ \ \ \ \ \ \ \ \ \ \ \ \ \ \ \ \
+l_3(\l,\m,\p)J_0+l_4(\l,\m,\p)J_{-1},
\end{eqnarray}
where $l_i(\l,\m,\p)$ for $1\leq i\leq 4$ are polynomials of $\l,\m,\p$. Using \eqref{LALALA}, by a little lengthy calculation, we also get
\begin{eqnarray}\!\!\!\!\!\!\!\!\!\!\!\!&\!\!\!\!\!\!\!\!\!\!\!\!\!\!\!&
\label{b-5}\ \ \ \ \ l_2(\l,\m,\p)=-\frac{1}{5}\l^2\big((15b-10a_0)\l^2+(51b-25a_0+15a_2)\l\m+(24b-10a_0)\l\p
\nonumber\\\!\!\!\!\!\!\!\!\!\!\!\!&\!\!\!\!\!\!\!\!\!\!\!\!\!\!\!&
\ \ \ \ \ \ \ \ \ \ \ \ \ \ \ \ \ \ \ \ \ \
+(15b-15a_0+35a_2)\m\p-10a_2\p^2\big),
\\
\!\!\!\!\!\!\!\!\!\!\!\!&\!\!\!\!\!\!\!\!\!\!\!\!\!\!\!&
\ \ \ \ \
l_4(\l,\m,\p)=\frac{1}{5}\l^2\big(b^2\l^4+4b^2\l^3\m+(4b^2-5b_2)\l^2\m^2+(4b^2-5b_2)\l\m^3-b^2\l^3\p
\nonumber\\\!\!\!\!\!\!\!\!\!\!\!\!&\!\!\!\!\!\!\!\!\!\!\!\!\!\!\!&
\ \ \ \ \ \ \ \ \ \ \ \ \ \ \ \ \ \ \ \
-(2b^2+5ba_2)\l^2\m\p+(3b^2-5ba_0-20b_2)\l\m^2\p-(b^2+5ba_0+15b_2)\m^3\p
\nonumber\\\!\!\!\!\!\!\!\!\!\!\!\!&\!\!\!\!\!\!\!\!\!\!\!\!\!\!\!&
\ \ \ \ \ \ \ \ \ \ \ \ \ \ \ \ \ \ \ \
+(6b^2+10b_2)\l^2\p^2+(26b^2-10ba_2+25b_2)\l\m\p^2-(b^2+5ba_0+15b_2)\m^2\p^2
\nonumber\\\!\!\!\!\!\!\!\!\!\!\!\!&\!\!\!\!\!\!\!\!\!\!\!\!\!\!\!&
\ \ \ \ \ \ \ \ \ \ \ \ \ \ \ \ \ \ \ \
+(11b^2+10b_2)\l\p^3+(6b^2-5ba_2+15b_2)\m\p^3\big)\label{b-6}.
\end{eqnarray}
From $(\ref{b-4})$, we know that $5\m+2\l+2\p$ must be a common factor of $l_2(\l,\m,\p)$ and $l_4(\l,\m,\p)$. Therefore by $(\ref{b-5})$ and $(\ref{b-6})$, we can deduce that $b$ must be zero.

Now assume $c\neq 0$.
By comparing the coefficients of  $J_{-1}$ on the both sides of \eqref{b-1}, we obtain
$$b_0=0,\,\ \ b_1=\frac15{b^2}-\frac52 bc^2,\,\ \ b_2=-\frac25b^2,\,\ \ b_3=-\frac12bc^2.$$ Then by comparing  the coefficients of  $J_{-1}$ on the both sides of \eqref{b-4} (this time, instead of \eqref{b-6}, we need to assume that
$l_4(\l,\m,\p)=(5\m+2\l+2\p)\sum_{i=0}^5d_i\l^{5-i}\p^i$), we obtain an equation, then by comparing coefficients of $\p^6,\p^5,\,\p^4$ in this equation, we
can obtain $d_5=d_4=b_4=0$ and $bc=0$. Thus $b=0$.
\hfill$\Box$\vskip5pt

Now we are ready to prove the main result of this paper.\vskip7pt

\ni{\it Proof of Theorem \ref{theo}.} By $(\ref{ver-4})$, Lemma
$\ref{lemm10}$ and Lemma $\ref{lemm11}$, we have
\begin{eqnarray}\label{MAMAMAM}\nonumber
\!\!\!\!\!\!\!\!\!\!\!\!\!\!\!\!\!\!\!\!\!\!\!\!\!\!\!\!&\!\!\!\!\!\!\!\!\!\!\!\!\!\!\!& \ \ \ \ \ \
\ \ [J{_{-1}}_{\l}J_{-1}]=0,\ \ \ \ \ \ \ \ \ \ \ \ \ \ \ \ \ \ \
\ \ \ \ \ \ \ \ \ \ \ \ \ \ \ \ \ [J{_{-1}}_{\l}J_0]=\l J_{-1},
\nonumber\\\!\!\!\!\!\!\!\!\!\!\!\!\!\!\!\!\!\!\!\!\!\!\!\!\!\!\!\!&\!\!\!\!\!\!\!\!\!\!\!\!\!\!\!&
\ \ \ \ \ \ \ \ [J{_0}_\l J_0]=(2\l+\p)J_0,\ \ \ \ \ \ \ \ \ \ \ \
\ \ \ \ \ \ \ \ \ \ \ \ \ \ \ [J{_{-1}}_{\l}J_1]=2\l
J_{0}+c\l^2J_{-1},
\nonumber\\\!\!\!\!\!\!\!\!\!\!\!\!\!\!\!\!\!\!\!\!\!\!\!\!\!\!\!\!&\!\!\!\!\!\!\!\!\!\!\!\!\!\!\!&
\ \ \ \ \ \ \ \ [J{_0}_{\l}J_1]=(3\l+\p)J_1+c\l^2J_{0},\ \ \ \ \ \
\ \ \ \ \ \ \ \ \ \ [J{_{-1}}_{\l}J_2]=3\l
J_1+3c\l^2J_{0}+c^2\l^3J_{-1},
\nonumber\\\!\!\!\!\!\!\!\!\!\!\!\!\!\!\!\!\!\!\!\!\!\!\!\!\!\!\!\!&\!\!\!\!\!\!\!\!\!\!\!\!\!\!\!&
\ \ \ \ \ \ \ \ [J{_{0}}_{\l}J_2]=(4\l+\p)
J_2+3c\l^2J_{1}+c^2\l^3J_{0},\ \ \ [J{_{1}}_{\l}J_1]=2(2\l+\p)
J_2-c(2\l\p+\p^2)J_{1}, \end{eqnarray} for some $c\in
\C.$ If $c\neq 0$, by replacing $J_{-1}$ and $J_1$ by $J'_{-1}=-c J_{-1}$ and $J'_1=-c^{-1}J_1$ respectively, we
can suppose $c=-1$. Thus, we can suppose $c$ is either equal to $0$
or $-1$. We want to prove
\begin{equation}\label{gc1-bracket-Final}
[J_{{m}_{\l}} J_{n}]= \left\{
\begin{array}{ll}
\Big((m+n+2)\l+(m+1)\p\Big)J_{m+n}&\mbox{if \ }c=0,\\[12pt]
\sum\limits_{s=0}^{m}\binom{m+1}{s+1}
(\l+\p)^{s+1}J_{m+n-s}-\sum\limits_{s=0}^{n}\binom{n+1}{s+1}
(-\l)^{s+1}J_{m+n-s}&\mbox{if \ }c=-1, \end{array}\right.
\end{equation}for $m,n\geq -1$. By \eqref{MAMAMAM}, we see
\eqref{gc1-bracket-Final} holds for all $m,n$ with ${\rm
max}\{m,n,n+m\}\le2$. Now inductively assume that for $N\ge2$,
\eqref{gc1-bracket-Final} holds for all $m,n$ with ${\rm
max}\{m,n,n+m\}\le N$. Denote the right hand side of
\eqref{gc1-bracket-Final} by $\sum_{k=-1}^{m+n}F_{m,n,k}(\l,\p)J_k$
for some polynomials $F_{m,n,k}(\l,\p)$. \vs{-10pt}Assume
\begin{eqnarray}\label{FFFFF-1}
&\!\!\!\!\!\!\!\!\!\!\!\!&[J_{{1}_{\l}}
J_{N}]-\mbox{$\sum\limits_{k=-1}^{N+1}$}F_{1,N,k}(\l,\p)J_k=\mbox{$\sum\limits_{k=-1}^N$}p_k(\l,\p)J_k,\\
\label{FFFFF-2} &\!\!\!\!\!\!\!\!\!\!\!\!&[J_{{0}_{\l}}
J_{N+1}]-\mbox{$\sum\limits_{k=-1}^{N+1}$}F_{0,N+1,k}(\l,\p)J_k=\mbox{$\sum\limits_{k=-1}^N$}
q_k(\l,\p)J_k,\\
\label{FFFFF-3} &\!\!\!\!\!\!\!\!\!\!\!\!& [J_{{-1}_{\l}}
J_{N+1}]-\mbox{$\sum\limits_{k=-1}^{N}$}F_{-1,N+1,k}(\l,\p)J_k=\mbox{$\sum\limits_{k=-1}^{N-1}$}
r_k(\l,\p)J_k\vs{-10pt},
\end{eqnarray}
for some polynomials $p_k(\l,\p),q_k(\l,\p),r_k(\l,\p)$. Applying the $\mu$-brackets
$J_{{0}_{\,\sc\m}},\,J_{{-1}_{\,\sc\m}}$ to \eqref{FFFFF-1} respectively,
using inductive assumption, we \vs{-10pt}obtain
\begin{eqnarray}\label{FFFFF-4}
\nonumber&\!\!\!\!\!\!\!\!\!\!\!\!&
-F_{1,N,N+1}(\l,\p)\mbox{$\sum\limits_{k=-1}^N$}q_k(\mu,\p)J_k=\mbox{$\sum\limits_{k=-1}^N$}
p_k(\l,\p)\mbox{$\sum\limits_{k'=-1}^k$}
F_{0,k,k'}(\mu,\p)J_{k'}\\
&\!\!\!\!\!\!\!\!\!\!\!\!&\phantom{-F_{1,N,N+1}(\l,\p)\mbox{$\sum\limits_{k=-1}^N$}q_k(\mu,\p)J_k}
=\mbox{$\sum\limits_{k=-1}^N\Big(\sum\limits_{k'=k}^N$}p_{k'}(\l,\p)F_{0,k',k}
(\mu,\p)\Big)J_k,
\\
\label{FFFFF-5} \nonumber&\!\!\!\!\!\!\!\!\!\!\!\!&
-F_{1,N,N+1}(\l,\p)\mbox{$\sum\limits_{k=-1}^{N-1}$}r_k(\mu,\p)J_k=\mbox{$\sum\limits_{k=-1}^{N}$}
p_k(\l,\p)\mbox{$\sum\limits_{k'=-1}^{k-1}$}
F_{-1,k,k'}(\mu,\p)J_{k'}\\
&\!\!\!\!\!\!\!\!\!\!\!\!&\phantom{-F_{1,N,N+1}(\l,\p)\mbox{$\sum\limits_{k=-1}^N$}q_k(\mu,\p)J_k}
=\mbox{$\sum\limits_{k=-1}^{N-1}\Big(\sum\limits_{k'=k+1}^N$}p_{k'}(\l,\p)F_{-1,k',k}
(\mu,\p)\Big)J_k\vs{-5pt}.
\end{eqnarray}
Note from the right hand side of \eqref{gc1-bracket-Final} that
$F_{1,N,N+1}(\l,\p)=(N+3)\l+2\p$, and $F_{-1,k',k}(\mu,\p)=\mu
F'_{k',k}(\mu)$ for some polynomial $F'_{k',k}(\mu)$ on $\mu$.
Comparing the coefficients of $J_k$ for $k=N,N-1,...,$ in
\eqref{FFFFF-5} shows that $\mu$ and $(N+3)\l+2\p$ are factors of $r_k(\mu,\p)$ and
$p_{k}(\l,\p)$ respectively. Thus $r_k(\mu,\p)=\mu r'_k(\mu,\p)$
and $p_k(\l,\p)=((N+3)\l+2\p)p'_k(\l,\p)$ for some $r'_k(\mu,\p)$
and $p'_k(\l,\p)$. Furthermore, we see that $r'_k(\mu,\p)=r'_k(\p)$
and $p'_k(\l,\p)=p'_k(\p)$ do not depend on $\l,\mu$. Thus in
\eqref{FFFFF-4}, if we replace $J_{N+1}$ by
$J_{N+1}-\sum_{k=-1}^Np'_k(\p)J_k$, we see that the right hand side
of \eqref{FFFFF-4} becomes zero, i.e., by re-choosing the generator
$J_{N+1}$, we can suppose all $p_k(\l,\p)=0$. Then \eqref{FFFFF-4}
and \eqref{FFFFF-5} show that all $q_k(\mu,\p),r_k(\mu,\p)$ are
zero. Hence \eqref{gc1-bracket-Final} holds for all $m,n$ with ${\rm
max}\{m,n,n+m\}=N+1$ and $m\le1$. To prove \eqref{gc1-bracket-Final}
for
 $2\le m\le n$ with
$m+n=N+1$, we use $J_m=\frac1{F_{1,m-1,m}(\mu,\p)}\big([J_{{1}_{\mu}}J_{m-1}]-\sum_{k=-1}^{m-1}F_{1,m-1,k}(\mu,\p)J_k\big)$,
and  Jacobi identity and induction on $m$. This proves
\eqref{gc1-bracket-Final} and Theorem \ref{theo}(1).
%
\hfill$\Box$ \vskip12pt

\cl{\bf\S5. \  Proof of Theorem \ref{theo}(2)
}\setcounter{section}{5}\setcounter{equation}{0}\setcounter{theo}{0} \vs{5pt}
\def\a{\alpha}\def\D{\Delta}

Assume $V$ is
a finitely freely $\C[\ptl]$-generated nontrivial  ${\rm gr\,}gc_1$-module. Regarding $V$ as a module over $\Vir$, by \cite[Theorem 3.2(1)]{CK}, we can choose a composition \vs{-5pt}series,
$$
V=V_{N}\supset V_{N-1}\supset\cdots\supset V_{1}\supset V_{0}=0\vs{-5pt},
$$
such that for each $i=1,2,...,N$, the composition factor $\overline{V}{}_{i}=V_{i}/V_{i-1}$ is either a rank one free module $M_{\Delta_i,\alpha_i}$ with $\D_i\ne0$, or else a $1$-dimensional trivial module $\C_{\a_i}$ with trivial $\l$-action and with $\ptl$ acting as the  scalar $\a_i$. Denote by $\bar v_i$ a $\C[\ptl]$-generator of $\overline{V}{}_{i}$ and $v_i\in V_i$ the preimage of $\bar v_i$.
Then $\{v_i\,|\,1\le i\le N\}$ is a $\C[\ptl]$-generating set of $V$, such that the $\l$-action of $J_0$ on $v_i$ is a $\C[\l,\ptl]$-combination of $v_1,...,v_i$.
\begin{lemm}\label{J-1-l-ac}For all $i\gg0$,  the $\l$-action of $J_i$ on $v_1$ is trivial, namely, ${J_i}_{\l}v_1=0$.
\end{lemm}
\noindent{\it Proof.~}~Assume  $i\gg0$ is fixed and suppose ${J_i}_{\l}v_1\ne0$, and let $k_i\ge1$ be the largest integer such that ${J_i}_{\l}v_1\not\subset V_{k_i-1}$. We consider the following possibilities.

\noindent{\bf Case 1.} $V_1=M_{\D_1,\a_1},\,\overline V_{k_i}=M_{\D_{k_i},\a_{k_i}}$.

We can \vs{-5pt}write
\begin{equation}\label{SSS}
{J_i}_\l v_1\equiv p_i(\l,\p)v_{k_i}\ ({\rm mod\,}V_{{k_i}-1})\mbox{ \ for some \ }p_i(\l,\p)\in\C[\l,\p]\vs{-5pt}.
\end{equation}
Applying the operator ${J_0}_\mu$ to \eqref{SSS}, we \vs{-5pt}obtain
\begin{equation}\label{SSS1}
p_i(\l,\mu\!+\!\p)(\a_{k_i}\!+\!\ptl\!+\!\D_{k_i}\mu)
\!=\!((1\!+\!i)\mu\!-\!\l)p_i(\l+\mu,\p)\!+\!(\a_1\!+\!\l\!+\!\p\!+\!\D_1\mu)p_i(\l,\p)\vs{-5pt}.
\end{equation}
Letting $\p=0$, we \vs{-10pt}obtain
\begin{equation}\label{==SSS1-PP}
p_i(\l,\mu)=-\frac1{\a_{k_i}+\D_{k_i}\mu}
\Big(((1+i)\mu-\l)p_i(\l+\mu,0)+(\a_1+\l+\D_1\mu)p_i(\l,0)\Big)\vs{-10pt}.
\end{equation}
Using this in \eqref{SSS1} with $\l=(1+i)\mu$ and $\p=-\a_{k_i}-\D_{k_i}\mu$, we \vs{-5pt}obtain
$$\begin{array}{ll}
(i+1)\big((\D_{k_i}{\sc}+1)\mu+\a_{k_i}\big)p_i\big((i+1-\D_{k_i})\mu-\a_{k_i},0\big)\\[5pt]
=\big((i+1-\D_1\D_{k_i})\mu+\a_1-\a_{k_i}\D_1\big)p_i\big((i+1)\mu,0\big)\vs{-5pt}.
\end{array}$$
Suppose $p_i(\l,0)$ has degree $m_i$. Comparing the coefficients of $\mu^{m_i+1}$ in the above equation, we obtain (note that the following equation does not depend on the coefficients of $p_i(\l,\mu)$\vs{-5pt})
\begin{equation}\label{==SSS1-PP===}
(i\!+\!1)(\D_{k_i}\!+\!1)(i\!+\!1\!-\!\D_{k_i})^{m_i}
\!=\!(i\!+\!1\!-\!\D_1\D_{k_i})(i\!+\!1)^{m_i}\vs{-5pt}.
\end{equation}
When $i$ is sufficient large, one can easily see that \eqref{==SSS1-PP===} cannot hold if $m_i>1$ (note that $\D_1,\D_{k_i}\ne0$, and $\D_{k_i}$ has only a finite possible choices since $1\le k_i\le N$). Thus $m_i\le1$ if $i\gg0$.
 Then from \eqref{==SSS1-PP}, we obtain that $p_i(\l,\mu)$ is a polynomial of degree $\le1$. Thus suppose $p_i(\l,\mu)=a_{i,0}+a_{i,1}\l+a_{i,2}\mu$. Then \eqref{SSS1} immediately gives $p_i(\l,\mu)=0$.\vs{5pt}

\noindent{\bf Case 2.} $V_1=\C_{\a_1},\,\overline V_{k_i}=M_{\D_{k_i},\a_{k_i}}$.

In this case, we can still assume \eqref{SSS}. Applying the operator ${J_0}_\mu$ to \eqref{SSS}, we \vs{-5pt}obtain
\begin{equation}\label{SSS2}
p_i(\l,\mu+\p)(\a_{k_i}+\ptl+\D_{k_i}\mu)=((1+i)\mu-\l)p_i(\l+\mu,\p)\vs{-5pt}.
\end{equation}
Letting $\mu=\p=0$, we obtain $p_i(\l,0)=0$. Then letting $\p=0$, we obtain $p_i(\l,\mu)=0$.

\noindent{\bf Case 3.} $V_1=M_{\D_1,\a_1},\,\overline V_{k_i}=\C_{\a_{k_i}}$.

In this case, since   $\p$ acts on $\bar v_{k_i}$ as the scalar $\a_{k_i}$, i.e., $\p v_{k_i}\equiv \a_{k_i} v_{k_i}\,({\rm mod\,}V_{{k_i}-1})$, we can writ\vs{-5pt}e
\begin{equation}\label{SSS3}
{J_1}_\l v_1\equiv p_i(\l)v_{k_i}\ ({\rm mod\,}V_{{k_i}-1})\mbox{ \ for some \ }p_i(\l)\in\C[\l].
\end{equation}
Applying the operator ${J_0}_\mu$ to \eqref{SSS3}, we \vs{-5pt}obtain
\begin{equation}\label{SSS3-1}
0=((1+i)\mu-\l)p_i(\l+\mu)+(\a_1+\l+\p+\D_1\mu)p_i(\l)\vs{-5pt}.
\end{equation}
By comparing the coefficients of $\p$, we immediately obtain $p_i(\l)=0.$

\noindent{\bf Case 4.} $V_1=\C_{\a_1},\,\overline V_{k_i}=\C_{\a_{k_i}}$.

As above, we immediately obtain $p_i(\l)=0$.
\hfill$\Box$\vs{5pt}

By induction on $j\le N$, we obtain ${J_i}_\l v_j=0$, i.e., the $\l$-action of $J_i$ is trivial. From this, we immediately obtain that the $\l$-action of ${\rm gr\,}gc_1$ on $V$ is trivial since ${\rm gr\,}gc_1$ is a simple conformal algebra. This proves Theorem
\ref{theo}(2).
\vs{10pt}
\small
\parskip=-1pt\baselineskip=3pt\lineskip=3pt

\noindent{\bf Acknowledgment.~}~We wish to thank the referee for suggesting the important improvements
 of the proof of Theorem \ref{theo}(1) by using
the $2$-cohomology (from which we obtain
Theorem \ref{Prop1}), and for suggesting the problem on whether ${\rm gr\,} gc_1$ has a
faithful representation on a finite $\C[\ptl]$-module (from which we obtain
Theorem \ref{theo}(2)).

\end{CJK*}
\end{document}